\documentclass[preprint]{imsart}
\RequirePackage[OT1]{fontenc}
\RequirePackage{amsthm,amsmath,natbib}
\RequirePackage{hyperref}
\RequirePackage{hypernat}
\setcitestyle{numbers}

\startlocaldefs
\numberwithin{equation}{section}
\theoremstyle{plain}
\newtheorem{Lemma}{Lemma}[section]
\newtheorem{thm}{Theorem}[section]
\newtheorem{Corollary}{Corollary}[section]
\newtheorem{Definition}{Definition}
\usepackage{epsfig}
\usepackage{epsf}
\usepackage{here}
\usepackage{color}
\usepackage{pifont}

\newcommand{\FF}{{{\cal F}}}
\newcommand{\Ob}{{{\cal O}}}
\newcommand{\X}{{{\cal X}}}

\newcommand{\Y}{{{\cal Y}}}
\newcommand{\Z}{{{\cal Z}}}
\newcommand{\R}{{{\cal R}}}
\newcommand{\N}{{{\cal N}}}
\newcommand{\W}{{{\cal W}}}
\newcommand{\LL}{{\cal L}}
\newcommand{\G}{{\cal G}}

\newcommand{\LRX}{{\cal L}^{\psi/\psi_0}_{{\cal R}|{\cal X}}}
\newcommand{\LOX}{{\cal L}^{\psi/\psi_0}_{{\cal O}|{\cal X}}}

\newcommand{\LX}{{\cal L}^{\theta/\theta_0}_{{\cal X}}}
\newcommand{\EZZ}{{\rm E}_{(\theta_0,\psi_0)}}

\newcommand{\EZZZ}{{\rm E}_{(\theta,\psi)}}
\newtheorem{Example}{Example}
\endlocaldefs

\begin{document}
\begin{frontmatter}

\title{ Likelihood inference for incompletely observed stochastic processes: ignorability conditions}
\runtitle{Ignorability conditions for stochastic processes}
\begin{aug}
\author{\fnms{Daniel} \snm{Commenges}\ead[label=e1]{daniel.commenges@isped.u-bordeaux2.fr}
 \and
\fnms{Anne} \snm{G\'egout-Petit}\ead[label=e2]{anne.petit@u-bordeaux2.fr}}

\address{Daniel Commenges, Epidemiology and Biostatistics Research Center, INSERM, Universit\'e Victor Segalen Bordeaux 2, 146 rue L\'eo Saignat, Bordeaux, 33076, France\\ Anne G\'egout-Petit, IMB. INRIA project CQFD,
Universit\'e Victor Segalen Bordeaux 2, 3, Place de la Victoire, 33076 BORDEAUX Cedex\\
\printead{e1,e2}} \runauthor{D Commenges and A G\'egout-Petit}
\affiliation{INSERM and Bordeaux University}
\end{aug}

https://sites.google.com/site/danielcommenges/

\begin{abstract} We develop a study of ignorability and conditions thereof for likelihood
inference in the framework of stochastic processes. We define a coarsening model for processes which  includes discrete-time observations as well as censored continuous-time observations and applies to continuous state-space processes as well as counting processes.
For preparing the work we recall formulas for manipulating  marginal and conditional
likelihood ratios (which can apply to stochastic processes). Ignorability is defined in terms of local equality
of two likelihood ratios. We give static conditions of ignorability and then dynamical conditions which are more
interpretable. We illustrate the use of the dynamical conditions of ignorability in problems of censoring,
missing data and joint modelling.
\end{abstract}


\tableofcontents
\end{frontmatter}

\section{Introduction}

Incomplete data are very common in statistics: when the mechanism leading to incomplete data (m.l.i.d.) is fixed a relatively simple likelihood can be written in general. Often the m.l.i.d. can not be considered as fixed and the question arises whether it can still be ignored.
Rubin \cite{Rubin} introduced the concept of ignorability for the simplest case in which the observation is a sequence of random variables and some of them are missing; he gave conditions
under which inference based on the assumption of fixed m.l.i.d. was valid, even in the case when in fact it
was not fixed. He established a typology of cases of missing data and showed in particular that the m.l.i.d. was
ignorable for likelihood inference in the case of missing at random (MAR) observations.

 In the framework of survival analysis the most frequent cases of incomplete data is right-censoring \cite{Kaplan,Cox} and  interval censoring \cite{Peto}. Conditions under which the conventional likelihood for right-censored
survival data was valid have been studied \cite{Laga,Kalb}. Andersen et al.
\cite{ABGK} developed the concept of independent censoring in the counting process framework. Heitjan and Rubin
\cite{Heitjan} also proposed some less conventional incomplete data cases which they called ``coarsening''. This topic
was also studied by Jacobsen and Keiding \cite{Jacobsen}, Gill et al. \cite{Gill} and Nielsen \cite{Nielsen}. The problematic of non-ignorable m.l.i.d.
has prompted the development of joint models, in which the m.l.i.d. was included, for instance in a model proposed by Diggle and Kenward \cite{Diggle}; see \cite{Thie} for a recent example.

The aim of this paper is to study ignorability in the context of stochastic processes: these processes may be
counting processes but also continuous state-space processes, such as diffusion processes. For giving a rigorous
treatment of that topic, we will need to rely on basic probability tools. First we will speak in terms of
likelihood ratio which is defined as a Radon-Nikodym derivative: this enables to manipulate likelihood ratios
for the observation of stochastic processes (for a review see \cite{Barndorff}). Local
equalities of $\sigma$-fields and of random variables will play an important role in the very definition of
ignorability and in the proofs and we recall a probability result described in \cite{Kall}:
\begin{Lemma}\label{Kallenberg}
Let the  $\sigma$-fields $\FF$, ${\cal G} \subset {\cal A}$ and
functions $\xi,\eta \in L^1$ be such that $A \cap \FF =A \cap
{\cal G} $ and $\xi=\eta$ a.s. on some set $A \in \FF\cap {\cal
G}$. Then ${\rm E}[\xi |\FF]={\rm E}[\eta |{\cal G}]$ a.s. on $A$.
\end{Lemma}
 Also, results on the likelihood of point
processes due to Jacod \cite{Jacod} will play a key role in several
proofs.

 We begin in section 2 by recalling the general definition of
the likelihood ratio and of marginal and conditional likelihoods which are valid for stochastic processes; a set of useful formulae is given. In section
3 a coarsening model for stochastic processes is given: it is represented by a stochastic process $R=(R_t)$ which indicates at each time whether the process of interest, $X$, is observed or not; if the process of interest is multivariate, $R$ can be multivariate. Then in section 4 we present a formulation of the
incomplete observation problem based on $\sigma$-fields. If $R$ is random the likelihood should include it, which complicates the inference problem.  We give a definition of ignorability: essentially the mechanism leading to incomplete data is ignorable if the full likelihood is equal to the likelihood of what is observed of $X$ in a model where the response process $R$ is fixed and equal to its observed value.    In section 5
we give static conditions of ignorability: the main condition called CAR(TCMP) (Definition 3) is that the likelihood of $R$ given $X$ depends only on observed quantities. In section 6 we give a dynamical condition called CAR(DYN) (Definition 5)  which is more interpretable
and usable than the general ones in some contexts: it says that the law of the response process only depends on what has been observed up to time $t$; this is made rigorous by  expressing this in terms of the equality of the compensator of the $R$ process in a filtration including $X$ and in the observed filtration.   Finally section 7 illustrates the use of the theory in
survival models, longitudinal data and joint models. Moreover we briefly illustrates different points along the paper by an example drawn from AIDS epidemiology. Section 8 is a short conclusion.

\section{Full, marginal and conditional likelihood}
Consider a measurable space $(\Omega, {\cal F})$ and a family of
measures $\{P_{\theta}\}_{\theta \in \Theta}$ absolutely
continuous relatively to a dominant measure $P_{\theta_0}$. For
$\X$ a sub-$\sigma$-field of $\FF$ the likelihood ratio on $\X$
is defined by:
$${\cal L}_{\cal X}^{\theta/\theta_0}=\frac{dP_{\theta}}{
dP_{\theta_0}}_{| {\cal X}}\hspace{5mm} {\rm a.s.} $$ where
$\frac{dP_{\theta}}{ dP_{\theta_0}}_{|{\cal X}}$ is the
Radon-Nikodym derivative of $P_{\theta}$ relatively to
$P_{\theta_0}$. Recall that $\frac{dP_{\theta}}{ dP_{\theta_0}}_{|
{\cal X}}$ is the ${\cal X}$-measurable random variable such that
$  P_{\theta}(F)=\int_F \frac{dP_{\theta}}{ dP_{\theta_0}}_{|
{\cal X}}dP_{\theta_0}, F\in {\cal X}$. If $\X \subset \Y \subset
\FF $ we have the fundamental formula \cite{Williams}:
$${\cal L}^{\theta/\theta_0}_{\X}={\rm E}_{\theta_0} [{\cal
L}^{\theta/\theta_0}_{\cal Y}|\X] \hspace{5mm} {\rm a.s.}.$$
Note that because conditional expectations and likelihood ratios
are defined a.s., all the equalities involving them are to be
understood as a.s., even if this is not specified for sake of
notational simplicity.

When $\FF$ is generated by two random elements $X$ and $Y$ and denote by $\X$ and $\Y$ the $\sigma$-fields they
generate respectively; thus $\FF=\X \vee \Y$; we will note ${\cal L}^{\theta/\theta_0}_{\cal F}= {\cal
L}^{\theta/\theta_0}_{\X \vee \Y}={\cal L}^{\theta/\theta_0}_{\X,\Y}$. The likelihoods
$\LL^{\theta/\theta_0}_{\X}$ and $\LL^{\theta/\theta_0}_{\Y}$ are called marginal likelihoods, and are linked to
the full likelihood by the conditional expectations: $\LL^{\theta/\theta_0}_{\X}={\rm E} [{\cal
L}^{\theta/\theta_0}_{\X,\Y}|\X]$ and $\LL^{\theta/\theta_0}_{\Y}={\rm E} [{\cal
L}^{\theta/\theta_0}_{\X,\Y}|\Y]$, as derives from the fundamental formula.

Conditional likelihoods can also be defined (see \cite{Hoffmann}); for brevity we do not recall the
definition. The conditional likelihood ratio of $Y$ given $\X$ will be denoted $\LL^{\theta/\theta_0}_{\Y|\X}$.
The following properties will be used in this paper:
\begin{eqnarray*}
{\rm i)}& \LL^{\theta/\theta_0}_{\Y|\X} \mbox{ is }  \X \vee\Y-{\rm measurable}\\
{\rm ii)}& \forall A \in \Y, \;\; {\rm E}_{\theta} \left[I_{A} | \X\right]=
{\rm E}_{\theta_0} \left[I_{A}\LL^{\theta/\theta_0}_{\Y|\X} | \X\right], a.s.\\
{\rm iii)}&  \LL^{\theta/\theta_0}_{\Y,\X}=\LL^{\theta/\theta_0}_{\Y|\X}\LL^{\theta/\theta_0}_{\X} \hspace{5mm}
{\rm a.s.}\\
{\rm iv)} &  \Y \subset \X \Longrightarrow \LL^{\theta/\theta_0}_{\Y,\X}=\LL^{\theta/\theta_0}_{\X} \iff
\LL^{\theta/\theta_0}_{\Y|\X}=1 \iff \LL^{\theta/\theta_0}_{\X|\Y}=
\frac{\LL^{\theta/\theta_0}_{\X}}{\LL^{\theta/\theta_0}_{_\Y}}\\
{\rm v)} &   \LL^{\theta/\theta_0}_{\X_1,\ldots,\X_m}=\LL^{\theta/\theta_0}_{\X_1}\prod_{k=2}^m\LL^{\theta/\theta_0}_{\X_k|\X_1,\ldots,\X_{k-1}}\\
{\rm vi)} & \X \vee \Y=\X \vee \Z \Longrightarrow \LL^{\theta/\theta_0}_{\Y|\X}=\LL^{\theta/\theta_0}_{\Z|\X}\\
\end{eqnarray*}
Note that ii) is the generalization of the main property of the likelihood ratio to conditional expectations;
thus $\LL^{\theta/\theta_0}_{\Y|\X}$ deserves its name of conditional likelihood ratio; it also implies
${\rm E}_{\theta_0} \left[\LL^{\theta/\theta_0}_{\Y|\X} | \X\right]=1$. Note also that under the assumption of a
family of equivalent measures, all the likelihoods are strictly positive a.s.

\section{A coarsening model for processes}
\subsection{Time coarsening: the $(R_t)$ process}
We define a  Time-Coarsening Model for Processes (TCMP) which is a general scheme where a process can be observed on a set of times, the set being possibly random. In many
real studies, in particular in epidemiology, we would have to
consider a sample of $n$ independent ``subjects'', to each of whom
a process $X^i=(X^i_t)=(X^i_t)_{t\ge 0}$ would be associated.
Since the likelihood would be the product of the individual
likelihoods, it is sufficient to consider only one process. We
first consider a process $X=(X_t)$ where $X_t$ takes values in
$\Re$, then we will extend the model to a multivariate process.
The main objective is to describe observation schemes for
processes in continuous time $t$, but $t$ may also be discrete so
that the results can be applied to finite collections of random
variables. We shall consider a response indicator process
$R=(R_t)$ taking value $1$ at $t$ if $X_t$ is observed and $0$
otherwise; this a generalization of the response indicator
variable introduced by Rubin \cite{Rubin}. This unifies different
concepts of censoring and observation of longitudinal data. Particular cases are:
\\ (i) right-censored survival data:
case where $X$ is a $0-1$ counting process and $R_t=1_{t\le C}$, where $C$ is a censoring variable;
\\ (ii) left-censored survival data:  case where $X$ is a $0-1$ counting
process and $R_t=1$ if $t\ge C$, $0$ otherwise;
\\ (iii) interval-censored survival data : case where $X$ is a $0-1$ counting
process  and $R_t=1$ if $t\in \{V_1, V_2, \ldots, V_m\}$, $0$ otherwise.
\\ (iv) Repeated measurements: case where $X$ has a continuous state
space and $R_t=1$ if $t\in \{V_1, V_2, \ldots, V_m\}$, $0$ otherwise.

\begin{figure}[!ht]
\begin{center}
\includegraphics{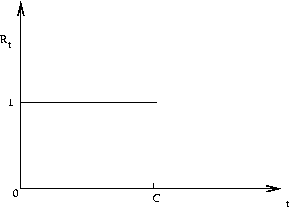}
\caption{right-censoring} \label{right-censoring}
\end{center}
\end{figure}

\begin{figure}[!ht]
\begin{center}
\includegraphics{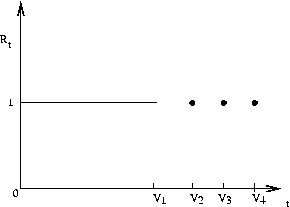}
\caption{Continuous monitoring followed by discrete-time visits} \label{non-standard}
\end{center}
\end{figure}

Note that $C$ in cases (i) and (ii), and $V_j, j=1,\ldots, m$ in (iii) and (iv) are
random variables. Cases (i) and (ii) illustrate a situation where $R$ is either right-
or left-continuous at each jump time and correspond to
observation in continuous time on some windows (see Figure 1);
case (iii) and (iv) correspond to observations in discrete time. In the
latter cases $R_t=1$ only on a finite or denumerable set of the
half line $[0,+\infty[$. The TCMP allows to represent a large
number of non-standard observation schemes. For instance the
subjects can be observed on windows separated by periods where no
observation are taken.  In most applications the process $X$ will
be either observed continuously on some periods or  observed only
at discrete time points, but the TCMP can represent a mixing of the
two types of observation for the same process: for instance a
subject could be observed in continuous time when he is at
hospital and at discrete times when he has left hospital (see
Figure 2). Also the process $X$ is not necessarily a 0-1 counting
process but may be for instance a more general counting process,
allowing recurrent events, or a process with continuous
state-space like a diffusion process.

\begin{Example} [Vital status, AIDS status]  The vital status, represented by a counting process, is in general observed in continuous time up to a censoring date (the end of the study). This is the conventional right-censoring case, case (i). The status with respect of a disease is often observed in discrete time. The AIDS status for an HIV infected person can be represented by a counting process $X$. Assume that the AIDS status is observed at visits at the hospital at times $V_1, V_2, \ldots,V_m$. This is a priori case (iii), that is, $R_t=1$ if $t\in \{V_1, V_2, \ldots, V_m\}$, $0$ otherwise. \end{Example}

That the TCMP includes censoring or coarsening models for random variables is obvious from the fact that to each
random variable $Y$ we can associate a counting process $X_t=1_{\{Y \ge t\}}$. We can define a
coarsening model for a random variable $Y$  by partitions ${\cal P}$ defined by intervals $A_j$ and
indicator $\nu_j$ which take values $1$ or $0$: if $\nu_j=1$, $Y$ is exactly observed on $\{Y\in A_j\}$, if $\nu_j=0$ it is
only observed that $Y$ falls in $A_j$. This is  equivalent to the TCMP for $X_t$ defined by
$R_t=\sum\nu_j1_{A_j}(t)$. In order to construct a random mechanism we can take a random partition, for instance determined by a set of random variables.

The model can be extended to the case of multivariate processes, as may be required by the observation of
several processes on the same ``subject''. So, we may consider that $X_t$ takes values in $\Re ^d, d \ge 1$ and
we may also consider a multidimensional response process $(R_t)$.

\begin{Example} [Joint vital and AIDS statuses] Consider a joint model of AIDS and vital status. The process $X=(X_1,X_2)$ is a bivariate counting process, $X_1$ representing AIDS and $X_2$ representing vital statuses respectively. The observation scheme can be represented by the bivariate response indicator process $R=(R_1,R_2)$, with $R_{1t}=1_{t\le C}$ and  $R_{2t}=1$ if $t\in \{V_1, V_2, \ldots, V_m\}$, 0 otherwise.\end{Example}

\noindent{\bf {\em Remark.}} The $(R_t)$ process is a generalization of the response indicator variables of Rubin \cite{Rubin}; van der Laan and Robins \cite{vdl} used similar response indicator processes for the case of counting processes in the right-censored and interval-censored cases. It is different from the filtering process proposed in Andersen et al. \cite{ABGK} (section III.4): for instance neither left-censoring nor discrete observation times can be treated with Andersen et al. filtering process. See also Arjas, Haara and Norros \cite{Arjas-H-N} and Arjas and Haara \cite{Arjas-H} for a detailed treatment of the filtering problem in the framework of marked point processes. The $(R_t)$ process can also be considered as a special case of the auxiliary random variable $G$ used by Heitjan and Rubin \cite{Heitjan} or Jacobsen and Keiding \cite{Jacobsen} if we interpret this random variable as a random element in a Skorohod space.

\subsection{Extension to vertical coarsening}
In addition to time-coarsening (which could be seen as ``horizontal coarsening''), ``vertical coarsening'' occurs if for $t$ such that $R_t=1$, $X_t$ is not completely observed, but a only coarsened version $Y_t=g(X_t,\varepsilon_t)$ where $\varepsilon_t$ is observed; here it is assumed for  simplicity that the processes $(\varepsilon_t)$ and $(X_t)$ are independent as well as the processes $(\varepsilon_t)$ and $(R_t)$.  Two cases, which are important in applications, are described: ``fixed vertical coarsening'' and noisy observation.

Fixed vertical coarsening is as follows: when $R_t=1$, $X_t$ is not completely observed, but according to a
fixed m.l.i.d.. This would apply to left-censored
observation of a biological marker due to a detection limit (such as HIV-RNA), as exemplified in section 7.2.3. In that case $(\varepsilon_t)=0$ and $Y_t=(1_{X_t> \eta}, \tilde X_t)$, where $\tilde X_t=\max (X_t,\eta)$.

It often occurs that only noisy observations are available. When $X_t$ is quantitative (i.d. can take any real value, case (iv) of section 3.1), a conventional model is that the observation that we make at  $t\in \{V_1, V_2, \ldots, V_m\}$ is noisy: typically we observe $Y(t)=X_t+\varepsilon_t$, with $\varepsilon_t$ representing a measurement error independent of all the other variables of the problem. This could also be applied in the case where $X_t$ is binary but there are classification errors.

\begin{Example} [CD4 counts]  \label{CD4} $X$ represents the concentration of T CD4 lymphocytes in an HIV infected person, and we have noisy observations of it $Y_j=X_{V_j}+\varepsilon_{j}$, $j=1,\ldots,m$. These noisy discrete-time observations are called the CD4 counts; see section \ref{randomobservationtimes}.\end{Example}

In the following the theory will be developed for the TCMP without vertical censoring; see Remark 3 of section 6 for the extension of the results to that case.

\section{Ignorability}
\subsection{The sigma-field representation for incomplete data}
A model for the random element $X$ is a family of measures
$\{P_{\theta}\}_{\theta \in \Theta}$ on a $\sigma$-field $\X$
generated by $X$ (for us, $X$ will be a stochastic process: $X_t$ takes
values in $\Re^d$ while the path of $X$ is an element of a Skorohod space). We will assume that the measures in the family
are equivalent and take $P_{\theta_0}$ as the reference measure.
If we observe $\X$ we will use the likelihood $\LL^{\theta
/\theta_0}_{\X}$ for inference about $\theta$. We will represent
the observed events by a $\sigma$-field ${\Ob}$. A general
definition of incomplete data is: $\X \not\subset \Ob$. A simple
case of incomplete data is when only a sub-$\sigma$-field of $\X$
has been observed: $\Ob \subset \X$.  A particular case occurs
when the m.l.i.d. can be represented by a TCMP with $R$ a
deterministic function. We shall denote by $\{R=r\}$  the event
$\{R_t=r_t, t\ge 0\}$ where $r_t$ is a particular path (an element
of the Skorohod space for instance). If $R$ is deterministic there
is a value (a path) $r$ such that $\{R=r\}=\Omega$.

\begin{Example}[Fixed right-censoring]
$R_t=1$ for  $t<c$ and $R_t=0$ for $t\ge c$, where $c$ is fixed. $R$ takes the fixed value $r$ such that $r_t=1_{t\le c}$.
In such a case we have  $\Ob=\X^r=\sigma(X_t, t\ge 0 :
r_t=1)=\sigma(r_tX_t,t\ge 0)$. In that case we have ${\cal
L}_{\Ob}={\cal L}_{\X^r}$ which is in general relatively easy to
compute. \end{Example}

\begin{Example} [Fixed right-censoring: the survival case] If $X$ is a counting process, Jacod's formula \cite{Jacod} can be used to obtain the likelihood; in the case where $X$ is a $0-1$ counting process with only one jump time $T$. The observed $\sigma$-field $\Ob$ can be seen as  generated by the variables $\delta=1_{\{T\le c\}}$ and $\tilde T= \min(T,c)$. The likelihood takes the form
$$\frac{\alpha_{\theta}^{\delta}(\tilde T)\exp[-\int_0^{\tilde T} \alpha_{\theta}(u) du ]}{\alpha_0^{\delta}(\tilde T)\exp[-\int_0^{\tilde T} \alpha_0(u) du ]},$$
 where $\alpha_{\theta}(.)$ (resp. $\alpha_0(.)$) are the risk functions under $P_{\theta}$ (resp. $P_0$).\end{Example}

Feigin \cite{Feigin} has given the likelihood for a diffusion process observed in continuous time.
If a process is observed at fixed times $v_1, \ldots, v_m$ the likelihood can be computed as the likelihood for observation of the vector of random variables $X_{v_1},\ldots, X_{v_m}$; if $X$ is a $0-1$ counting process this case has been denoted ``interval censoring'' and the likelihood is simple to compute (see \cite{Peto}, \cite{Alioum}); if $X$ is a Gaussian process, the vector  $X_{v_1},\ldots, X_{v_m}$ has a normal distribution which makes the likelihood easy to compute.

If $R$ is not fixed, the above definition of $\Ob$ is meaningless; we must include $R$ in the description of the
problem. We shall consider a larger $\sigma$-field $\FF=\X \vee \R$, where $\R$ is the $\sigma$-field generated by
$R$ for right- or left-continuous processes that is: $\R=\sigma(R_t, t\ge 0)$; if $R$ takes value $1$ at only a
finite (or denumerable) number of times (corresponding to the discrete observation case) we can take $\R$ as
generated by the counting process counting the number of observation times.
We consider that $R$ is observed (see Remark 2) so that a representation of $\Ob$ is $\Ob=\sigma(R_t X_t, R_t, t
\ge 0)$. In section 6 which develops a
dynamical approach to the problem, we shall define adequate filtrations; for instance if $R$ is c\`adl\`ag (right-continuous with left-hand limits) the observed filtration $(\Ob_t)$ will be the family of $\sigma$-fields $\Ob_t=\sigma(R_uX_u,R_u, 0\le u \le t)$. We have of course $\Ob \subset \FF$; generally we have  $\X \not\subset \Ob$ (incomplete data) and $\Ob
\ne \X^r $ (the observation is not a predetermined subset of values of $X$).

\noindent{\bf {\em Remark 1.}}  We might think that we could define an interesting $\sigma$-field by
$\sigma(R_t X_t, t \ge 0)$ which could take the role of the notation $x_{obs}$ used in most of the literature in
missing data (for instance \cite{Kenward}); however if $X_t \ne 0$ for all $t$, the latter
$\sigma$-field is equal to $\Ob$. When $R$ is random it is not possible to disentangle the observed part of $X$
from $R$; only the realized value of $R$ effects a partition between  $x_{obs}$ and $x_{mis}$. This is in fact
the meaning of $\X^r$ which is the observed part of $\X$ when $R=r$.

\noindent{\bf {\em Remark 2.}} It is natural to say that $R$ is observed: for each $t$ we know whether we observe $X_t$ or not. There is
however an important case where this natural assumption does not hold: in survival analysis, $X$ is a $0-1$
counting process, so after a jump has been observed there is no need for observation anymore; so we may ignore
whether we would have observed the process if it had been necessary. The simplest way to get out of this problem
is to put $R_t=1$ if $X_t$ is known. More generally assume that there is an absorbing state $a$ such that if
$X_t=a$, $X_{t+u}=a$ for $u>0$; define the $\Ob_t$-stopping time $T=\inf \{t: X_t=a \mbox{ and } R_t=1\}$. By convention put
$R_{T+u}=1, u>0$. This part of the law of $R$ is anyway unidentifiable. In the remaining of the paper we will
consider that $R$ is observed.

\noindent{\bf {\em Remark 3.}} If $X_t$ is multivariate and its components may be differently coarsened, then $R_t$ is multivariate and $R_tX_t$ must be interpreted as the scalar product of the two vectors.

\subsection{Model, notations and observed likelihood}
From now on, a  model for the random element $(X,R)$ is a family
of measures $\{P_{(\theta, \psi)}\}_{(\theta, \psi) \in \Theta
\times \Psi}$  on a measurable space $(\Omega, \FF)$. $X$ (resp.
$R$) takes values in a measurable space $(\Xi, \xi)$ (resp.
$(\Gamma, \rho)$). For us $X$ and $R$ will be  d-dimensional
c\`{a}dl\`{a}g stochastic processes, so $(\Xi, \xi)$ and  $(\Gamma, \rho)$ are Skorohod
spaces endowed with their Borel $\sigma$-fields.
 The parameter spaces $\Theta$
and $\Psi$ need not be finite dimensional. We will assume
that the measures in the family are equivalent and take
$P_{(\theta_0,\psi_0)}$ as the reference measure. $P_{\theta}$ is
the restriction of $P_{(\theta, \psi)}$ to $\X$: that is, the marginal probability of $X$ does not depend on $\psi$. The additional
parameter $\psi$ will be considered as a nuisance parameter. We
assume implicitly a "Non-Informativeness" assumption in the
coarsening mechanism, which is :

\begin{equation}\label{noninfo}
  P_{(\theta_1, \psi)}(A|\X)=P_{(\theta_2, \psi)}(A|\X)\mbox{, a.s., } A
\in \R \mbox{ for all } \theta_1,  \theta_2, \psi
\end{equation}

In words, the conditional probability of $R$ given $X$ does not depend
on $\theta$.
 This has an important consequence in terms of likelihood
ratio. The latter assumption can also be written as: ${\rm E}_{(\theta_1,
\psi)}(1_A|\X)={\rm E}_{(\theta_2, \psi)}(1_A|\X)$, $A \in \R$,
which remembering property ii) of the conditional likelihood is
equivalent to ${\rm E}_{(\theta_0, \psi_0)}(1_A{\cal
L}^{(\theta_1,\psi )/ (\theta_0,\psi_0)}_{\R|{\cal X}}|\X)={\rm
E}_{(\theta_0, \psi_0)}(1_A{\cal L}^{(\theta_2,\psi) /
(\theta_0,\psi_0)}_{\R|{\cal X}}|\X)$. This in turn implies that
${\cal L}^{(\theta_1,\psi) /( \theta_0,\psi_0)}_{\R|{\cal
X}}={\cal L}^{(\theta_2,\psi) /( \theta_0,\psi_0)}_{\R|{\cal X}}$
and we will denote this common value ${\cal L}^{\psi /
\psi_0}_{\R|{\cal X}}$. Moreover it can be proved (and is
intuitive) that ${\cal L}^{\psi_0 / \psi_0}_{\R|{\cal X}}=1$,
a.s..

From now on, we fix our notation for $\Ob$ and $\X^r$ by
$\Ob=\sigma(R_t X_t, R_t, t \ge 0)$ and  for $r=(r_t)$ a
deterministic path of $R$,  $\X^r=\sigma(X_t, t\ge 0 :
r_t=1)=\sigma(r_tX_t,t\ge 0)$. Note that (\ref{noninfo})  remains
true even if $A \in \Ob$. The inferential likelihood is ${\cal
L}^{(\theta,\psi) /(\theta_0,\psi_0) }_{\Ob}$ (which is
$\Ob$-measurable and thus can be computed from the observations)
and the fundamental property yields:
$${\cal L}^{(\theta,\psi) /(\theta_0,\psi_0) }_{\Ob}=\EZZ [{\cal L}^{(\theta,\psi) /(\theta_0,\psi_0) }_{\X, \R}|\Ob]=\EZZ [ \LRX \LX|\Ob]$$

\begin{Example} [Random right-censoring in the survival case] In the survival case we can distinguish the computation of the likelihood on the observed events $\{\delta=1\}$ and $\{\delta=0\}$. On $\{\delta=1\}$ $X$ is completely observed so that we can put $\LX$ out of the conditional expectation getting ${\cal L}^{(\theta,\psi) /(\theta_0,\psi_0) }_{\Ob}=\EZZ [ \LRX|\Ob] \LX$. Since $\LRX=\frac{P_{\psi}(C>T)}{P_{\psi_0}(C>T)}$ it is known on $\{\delta=1\}$ and can be put out of the conditional expectation. Moreover it does not depend on the parameter of interest $\theta$, so we could get rid of it for inference about $\theta$. On $\{\delta=0\}$ however, the conditional expectation is $\EZZ [ \LRX \frac{f_{\theta}(T)}{f_{0}(T)}|C,T>C]$, where $f_{\theta}(.)$ is the density function of $T$; since $X$ is not completely observed none of the terms can get out of it without additional assumptions. \end{Example}
\subsection{Definition of ignorability}

If the m.l.i.d. (represented by $R$) is random it may still be tempting to ignore it, treating it as fixed, and
use for inference  $\LL_{\X^r}$ which is relatively easy to compute and does not depend on $\psi$.

\begin{Definition}
In the TCMP, the likelihood ratio ignoring the m.l.i.d. is the likelihood ratio $\LL^{\theta /
\theta_0}_{\X^r}$ obtained under the assumption that $R$ is fixed at its observed value $r$.
\end{Definition}

When the fixed m.l.i.d. assumption does not hold, the question arises to know in which cases (if ever)
$\LL^{\theta / \theta_0}_{\X^r}$ leads to the same inference about $\theta$ as $\LL^{\theta,\psi_0/
\theta_0,\psi_0}_{\Ob}$, if the true value of $\psi$ was known to be $\psi_0$ or $\LL^{\theta,\psi/
\theta_0,\psi_0}_{\Ob}$ if $\psi$ had to be estimated.

For defining ignorability we face the problem  that both
Radon-Nikodym derivatives and conditional expectations are defined
almost surely; for some results we must restrict the theoretical
framework to measures giving a non null probability to a
denumerable set of trajectories of $R$.

{\bf Assumption A1} There exists a denumerable set $(r_1, r_2, \ldots)$ such that $P_{\theta}(R=r_i)>0$ and $\sum_i
P_{\theta}(R=r_i)=1$, for all $\theta$.

 This is a theoretical limitation but this has no impact on application since in practice
the times are always rounded. Some of the results below will need this restriction, other will not.

\begin{Definition}[Ignorability]
The m.l.i.d. will be called {\bf ignorable on {\em r} } if
$\LL^{(\theta,\psi/ (\theta_0,\psi_0)}_{\Ob} = U \LL^{\theta /
\theta_0}_{\X^r}$ \hspace{3mm} {a.s.} on $\{R=r\}$ for all
$\theta $, whatever $\psi$ and $\psi_0$, where $U$ is random variable not depending on $\theta$. It will be called {\bf ignorable} if
assumption A1 holds and
 the m.l.i.d. is ignorable on $r$ for all values of $r$.
\end{Definition}

\begin{Example} [Random right-censoring in the survival case] \label{lik-rightcensor} In the survival case on $\{\delta=0\}$ and if we observe $C=c$ the likelihood ignoring the m.l.i.d. is
$\LL^{\theta /\theta_0}_{\X^r}=S_{\theta}(\tilde T)/S_{\theta_0}(\tilde T)$. Thus, ignorability holds if $\EZZ [ \LRX \frac{f_{\theta}(T)}{f_{0}(T)}|C,T>C]=U\frac{S_{\theta}(\tilde T)}{S_{\theta_0}(\tilde T)}$, where $S_{\theta}(.)$ is the survival function of the distribution of $T$ under $P_{\theta}$.
\end{Example}

\noindent {\bf {\em Remark.}}
 If the m.l.i.d. is  ignorable it is then obvious that $\LL^{\theta
/ \theta_0}_{\X^r}$ and $\LL^{(\theta,\psi_0 )/ (\theta_0,\psi_0)}_{\Ob}$ lead to the same inference about
$\theta$.

\section{Static conditions of ignorability for the TCMP}

We give a first fact which does not seem to have been noted
previously in a general context.

\noindent {\bf{\em Fact.}} Ignorability on $\{R_t=1, t\ge
0\}$ always holds.

\noindent {\bf{\em Proof.}} On the
event $\{R_t=1, t\ge 0\}$, $X$ and $R$ are observed so that $\FF = \Ob$;
thus from Lemma \ref{Kallenberg} we have
  $\LL^{(\theta,\psi) /(\theta_0,\psi_0)
}_{\Ob}=\EZZ(\LL^{(\theta,\psi) /(\theta_0,\psi_0) }_{\FF}|\Ob)=\EZZ(\LL_{\FF}|\FF)=\LL^{(\theta,\psi)
/(\theta_0,\psi_0) }_{\FF}= \LRX \LX$ (the last equality comes from property iii) of the likelihood); for $\psi=\psi_0$ we retrieve the likelihood $\LX$ for the complete
observation of $X$ and ignorability holds on $\{R_t=1,t\ge 0\}$.\hfill \ding{113}

We shall now study ``static'' conditions of ignorability, in contrast with the ``dynamic'' conditions of the next section.
Gill \& al \cite{Gill} have introduced two conditions of ignorability: CAR(REL) (Relative Coarsening At Random) and CAR(ABS) (Absoluter Coarsening At Random); these conditions were further studied by Nielsen \cite{Nielsen}. We give two conditions of ignorability in
the TCMP framework. The first one is an adaption
of CAR(REL); the second one is stronger
but original and doesn't imply CAR(ABS).

\begin{Definition}[CAR(TCMP)]
We will say that CAR(TCMP) holds for the couple $(X,R)$ if $\LRX$
is $\Ob$-measurable for all $(\theta, \psi)$.
\end{Definition}

\begin{Example}[$R$ and $X$ independent]
If $R$ and $X$ are independent  $\LRX= {\cal L}^{\psi/\psi_0}_{{\cal R}}$; since $R$ is observed ${\cal L}^{\psi/\psi_0}_{{\cal R}}$ is $\Ob$-measurable.
\end{Example}

With our notations the condition CAR(REL) is "there exists a version of $\LOX(o,x)$ such that the
mapping $x \rightarrow \LOX(o,x)$ is constant for all $x$
compatible with $o$", where $o$ is an elementary event of $\Ob$, i.e. $o=(r,rx)$.
First note that with our assumptions, $x$ is compatible with $o=(r,y)$ if $rx=y$ and
 second note that $\LL_{\R|\X}(r ; x)=\LL_{\Ob|\X}(o ; x)$ if
 $o=(r,rx)$ (by property vi) and because $\R\vee \X=\Ob \vee \X$.
So in the TCMP setting the condition CAR(REL) becomes
"there exists a version of $\LRX(o,x)$ such that :
for all  $r$ and $ (x,x') $
verifying $rx=rx'$ then  $\LRX(r,x)= \LRX(r,x')"$.

\begin{thm}
CAR(TCMP) is equivalent to CAR(REL) in the TCMP setting.
\end{thm}

\noindent {\bf {\em Proof.}} If CAR(REL) is true then for all $(r,x)$, $\LRX(r ;
x)=\LRX(r ; rx)$ (by taking $x'=rx$ in the above formula) and then
$\LRX$ is a function
 of $(R,RX)$ and so $\LRX$ is $\Ob$-measurable that is CAR(TCMP).
 Conversely, if $\LRX$ is $\Ob$-measurable, there is a version of it
 which is constant on the atoms of $\Ob$ i.e. on the set of the
 form $O_{r,y}=\{(r,x) \mbox{ such that } x \mbox{ verifies } rx=y
 \}$ and we have CAR(REL).\hfill \ding{113}

The next theorem shows that CAR(TCMP) implies a factorization of
the likelihood in two parts : one which depends on $\psi$ and the
second one on $\theta$; this may be called ``weak ignorability''.

\begin{thm}[Factorization]
If the couple $(R,X)$ satisfies CAR(TCMP) then we have
$\LL^{(\theta, \psi) /(\theta_0, \psi_0)}_{\Ob}=\LRX
\EZZ(\LX|\Ob)$ and $\EZZ(\LX|\Ob)$ does not depend on $\psi_0$.
\end{thm}

\noindent {\bf {\em Proof.}} A proof can be obtained using the previous theorem
and the fact that CAR(REL) has been proved to imply a similar
factorization theorem. However a direct proof is quite simple in
our formalism. From the decomposition formula we have
$\LL^{(\theta, \psi) /(\theta_0, \psi_0)}_{\FF}=\LRX \LX$ and we
can pull $\LRX$ out of the conditional expectation using CAR(TCMP)
thus obtaining $\LL^{(\theta, \psi) /(\theta_0,
\psi_0)}_{\Ob}=\LRX \EZZ(\LX|\Ob)$. It only remains to prove that
the last term does not depend on $\psi_0$. We have
\begin{eqnarray*}\LL^{(\theta,\psi)/(\theta_0,\psi_0)}_{\Ob}&=&\LL^{(\theta,\psi)/(\theta_0,\psi_1)}_{\Ob}\LL^{(\theta_0,\psi_1)/(\theta_0,\psi_0)}_{\Ob}\\
&=& {\cal L}^{\psi/\psi_1}_{{\cal R}|{\cal X}} {\rm E}_{(\theta_0,\psi_1)}[{\cal L}^{\theta/\theta_0}_{{\cal X}}|\Ob]{\cal L}^{\psi_1/\psi_0}_{{\cal R}|{\cal X}} {\rm E}_{(\theta_0,\psi_0)}[{\cal L}^{\theta_0/\theta_0}_{{\cal X}}|\Ob]\\
&=& \LRX {\rm E}_{(\theta_0,\psi_1)}[{\cal L}^{\theta/\theta_0}_{{\cal X}}|\Ob]
\end{eqnarray*} \hfill \ding{113}

As a result, we can use $\EZZ(\LX|\Ob)$ for inference on $\theta$, getting rid of the nuisance parameter $\psi$.

\begin{Example} [Random right-censoring in the survival case] Continuing Example \ref{lik-rightcensor} we have on $\delta=0$, the likelihood if CAR(TCMP) holds can be written: $\EZZ [ \LRX \frac{f_{\theta}(T)}{f_{0}(T)}|C,T>C]=\LRX \EZZ [ \frac{f_{\theta}(T)}{f_{0}(T)}|C,T>C]=\LRX \frac{S_{\theta}(\tilde T)}{S_{\theta_0}(\tilde T)}$. Thus the observed likelihood is proportional to the likelihood ignoring the m.l.i.d.; it is strictly equal if we take $\psi=\psi_0$, our definition of ignorability.
\end{Example}

The following condition allows obtaining a more precise result locally.

\begin{Definition}[CAR(TCMP)-loc]\label{carloc}
We will call CAR(TCMP)-loc on $r$ the condition:  ${\cal
L}^{(\theta, \psi) /(\theta_0, \psi_0)}_{\R|{\cal X}}={\cal
L}^{(\theta, \psi) /(\theta_0, \psi_0)}_{\R|{\cal X}^r}$ a.s. on
$\{R=r\}$.

\end{Definition}

We have:
\begin{thm}\label{loc}
 CAR(TCMP)-loc on
$r$ implies ignorability on $r$.
\end{thm}

Details about this concept  and proof of Theorem \ref{loc}  are given in Appendix 1.

\section{Dynamical conditions of ignorability}

 Assume that we can represent $(R_t)$ by a marked point process $(N_t)$;
 if the component $R_h$ of $R$ takes value $1$ at isolated points $V_{h,1}, V_{h,2}, \ldots$
 (discrete-time observations), then its associated counting
 process $N^d_{ht}$ is  defined by
 $N^d_{h,t}=\sum_{j=1}^{\infty}1_{\{V_{h,j}\leq t\}}$. If the component
 $R_{h,t}$ is c\`{a}dl\`{a}g, it can be written as $R_{h,t}=\sum_{k=0}^{\infty}1_{[V_{h,2k},V_{h,2k+1}[}$
 and its associate counting process $N^c_{h,t}=\sum_{j=0}^{\infty}1_{\{V_{h,j}\leq t\}}$.
 If the component $R_h$ of $R$ is a mixing of the two kinds of
  observation,
  (discrete- and continuous-time) then at points of discontinuities we have
  to add a mark which indicates if the jump affects the discrete
  part of $R_h$ represented by $N^d_{h,t}$ or the c\`{a}dl\`{a}g one
  represented by $N^c_{h,t}$. In all cases $R$ can be represented by
  a marked point process $(N_t)$.
Denote by $(\N_t)$ the self-exciting filtration of $(N_t)$.  We
define the filtration $(\Ob_t)$ as the family of $\sigma$-fields
$\Ob_t=\sigma(N_u,R_uX_u, 0\le u \le t)$. Let us call
$\Lambda^{\Ob,N}=(\Lambda^{\Ob,N}_t)$  and
$\Lambda^{\FF^*,N}=(\Lambda^{\FF^*,N}_t)$ the compensators of $N$
in the filtrations $(\Ob_t)$ and $(\FF^*_t)$ respectively and for
probability $P_{\theta,\psi}$, where $(\FF^*_t)$ is the family of
$\sigma$-fields $\FF^*_t= \X \vee \Ob_t, t\ge 0$. The compensators
generally depend on $(\theta, \psi)$ but we omit this for
notational simplicity; the compensators for $P_{\theta_0,\psi_0}$
will be denoted $\Lambda_0$. In the following we will assume that
there exists a fixed time $\tau$ such that $R_{\tau+u}=R_{\tau}$
and  that there is no explosion of the process on $[0,\tau]$ (i.e.
$\sum_{t \leq \tau}\Delta R_t <\infty$) (thus
$\Lambda^{N}_{\tau+u}=\Lambda^{N}_{\tau}, u\ge 0$, for any
filtration).

\begin{Definition}[CAR(DYN)]

We will denote CAR(DYN) the condition:

$\rm{ \hspace{10mm}} \forall \;(\theta, \psi)\;\;  $ we have under
$P_{(\theta, \psi)}: \;\; (\Lambda^{\Ob,N}_t) =
(\Lambda^{\FF^*,N}_t)$, (up to indistinguishability).
\end{Definition}

\noindent{\bf {\em Remark 1.}} This is an absolute condition in the sense of Gill et al. \cite{Gill} and Nielsen \cite{Nielsen} since the CAR(DYN) criterion is defined for each probability separately while the definition of CAR(TCMP) bears on Radon-Nikodym derivatives.

\noindent{\bf {\em Remark 2.}} For this condition to lead really to ignorability the model for the mechanism leading to missing data must be well specified. That is, the condition $(\Lambda^{\Ob,N}_t) = (\Lambda^{\FF^*,N}_t)$ must hold under the true probability.

\begin{Example}[Discrete visit times] Suppose that $N$ represents visit times. This process has the compensators $(\Lambda^{\Ob,N}_t) $ and  $(\Lambda^{\FF^*,N}_t)$ in the filtrations $(\Ob_t)$ and $(\FF^*_t)$ respectively. CAR(DYN) holds if these compensators are equal. In words, the next visit time depends only on what has been observed at previous visit times. For instance if we are studying the concentration of CD4 T lymphocytes, this is the case if the next visit time is fixed as a function of the observed CD4 counts (see section 7.2.2). \end{Example}

It appears that CAR(DYN) is stronger than CAR(TCMP).

\begin{thm}
CAR(DYN) implies CAR(TCMP).
\end{thm}

\noindent {\bf {\em Proof.}} Let us write the likelihood for $N$ and $X$. The
filtration $(\FF^*_t)$ is the self-generated filtration for $N$
when $\FF^*_0=\X$. Thus we have using Jacod's formula $
\LL_{\FF^*}=\LX\phi(\Lambda^{\FF^*,N}_u,\Lambda^{\FF^*,N}_{0u},
N_u, u\ge 0)$, where $$\phi(\Lambda, \Lambda_0,
N)=\phi(\Lambda_u,\Lambda_{0u}, N_u , 0\le u)=\prod_{t \ge 0}\prod
_h \Biggl(\frac{d  \Lambda _h}{d\Lambda_{0 h}}_t\Biggr)^{\Delta
N_{ht}} \frac {\prod_{t \ge 0: \Delta N_{.t }\ne 1} (1-d  \Lambda
_{.t})}{\prod_{t \ge0: \Delta N_{.t} \ne 1} (1-d  \Lambda_{0
.t})},$$ where $N_.=\sum _h N_h$ and $\Lambda_.= \sum _h
\Lambda_h$. Identifying with the decomposition $\LL_{\FF^*}=\LX
\LRX$ we find that $\LRX=\phi(\Lambda^{\FF^*,N},
\Lambda^{\FF^*,N}_0, N)$. With CAR(DYN) this is equal to $\phi
(\Lambda^{\Ob,N}, \Lambda^{\Ob,N}_0, N)$ which is
$\Ob$-measurable, and thus CAR(TCMP) holds.\hfill \ding{113}

\begin{thm} \label{Dyn>loc}
CAR(DYN) implies ignorability.
\end{thm}

\noindent {\bf {\em Proof.}} We have to prove that CAR(DYN) implies CAR(TCMP)-loc for all $r$; then using Theorem \ref{loc} we will have that CAR(DYN) implies ignorability. Let us define $\FF^r_t=\X^r \vee \N_t$ and
$\FF^r=\X^r \vee \N$ then for all $t \in [0,\tau]$, $\FF^r_t
\subseteq \FF^*_t$. As seen in the previous proof, ${\cal
L}^{(\theta, \psi) /(\theta_0, \psi_0)}_{\R|{\cal
X}}=\phi(\Lambda^{\FF^*,N}, \Lambda^{\FF^*,N}_0, N)$ and ${\cal
L}^{(\theta, \psi) /(\theta_0, \psi_0)}_{\R|{\cal X}^r}=
\phi(\Lambda^{\FF^r,N}, \Lambda^{\FF^r,N}_0, N)$. So the proof
will be complete if $\forall t \in [0,t]$,
$\Lambda_t^{\FF^r,N}1_{\{R=r\}}= \Lambda_t^{\FF^*,N}1_{\{R=r\}}$.
We get
 \begin{eqnarray}
 \Lambda_t^{\FF^r,N} 1_{\{R=r\}} & = & E[\Lambda_t^{\FF^*,N} |\FF^r_t]1_{\{R=r\}} \label{ino}\\
 & = & E[\Lambda_t^{\Ob,N} 1_{\{R_{. \wedge t} =r_{. \wedge
t}\}}|\FF^r_t] 1_{\{R=r\}} \label{dyn}\\
& = & \Lambda_t^{\Ob,N}1_{\{R=r\}} \label{mesu}\end{eqnarray}

 (\ref{ino}) is due to the innovation theorem, (\ref{dyn}) is due
 to CAR(DYN) and the fact that
$\{R_{. \wedge t} =r_{. \wedge t}\} \subset \{R=r\}$ and $\{R_{.
\wedge t} =r_{. \wedge t}\}\in \Ob_t \cap \FF^r_t$. At last the
local inclusion $\{R_{. \wedge t} =r_{. \wedge t}\}\cap \Ob_t
\subseteq \{R_{. \wedge t} =r_{. \wedge t}\}\cap \FF^r_t$ implies
the $\FF^r_t$-measurability of $\Lambda_t^{\Ob,N} 1_{\{R_{. \wedge
t} =r_{. \wedge t}\}}$ and (\ref{mesu}). Using again CAR(DYN), we
have the desired equality.   \hfill \ding{113}

\begin{Corollary} \label{Rprev}
If $R$ is $\Ob_t$-predictable, CAR(DYN) holds and thus ignorability holds.
\end{Corollary}

\noindent {\bf {\em Proof.}} The associated counting process $N$ is itself $\FF^*_t$-predictable. Thus the Doob-Meyer decomposition is $N_t=N_t+0$ in both $(\Ob_t)$ and $(\FF^*_t)$. It follows that $\Lambda^{\FF^*,N}=\Lambda^{\Ob,N}=N$ which is CAR(DYN). By Theorem \ref{Dyn>loc} and \ref{loc} ignorability holds.\hfill \ding{113}

In the framework of counting processes, Andersen et al. \cite{ABGK} have proposed a criterion of independent
right-censoring. We adapt their criterion  to right-continuous censoring processes and we restrict to the univariate case for simplicity. In that case $R_t=1_{\{t<C\}}$
where $C$ is a censoring variable and $N_t=1-R_t$. Let $X$ be a counting process and $\Lambda^{\X,X}$ its
compensator in the self-generated filtration $(\X_t)$. Let $(\FF_t)$ the filtration generated by both $X$ and $N$.
We have independent right-censoring if the compensator of $X$ is the same in the filtration including
information on the censoring that is: $\Lambda^{\FF,X}=\Lambda^{\X,X}$. In the context where both the independent censoring and CAR(DYN) apply, the following theorem says that they are equivalent.

\begin{thm}[Equivalence of independent censoring and CAR(DYN)]\label{idpcensor}
Let $X$ be a counting process which admits a c\`ag (left-continuous) intensity $\lambda^{\X,X}$ in the self-generated filtration $(\X_t)$.
Consider a right-continuous right-censoring process of $X$ satisfying  CAR(DYN); then this is an  independent censoring. Inversely, independent censoring implies CAR(DYN).
\end{thm}
The proof is given in Appendix 2.

In some situations it is natural to consider response indicator
processes which are left-continuous; case II right-censoring is an
example. This is developed in Appendix 3.

\noindent{\bf {\em Remark 3.}} All the theoretical results can be extended to TCMP with vertical coarsening (section 3.2) by by replacing $X_t$ by $Y_t=g(X_t,\varepsilon_t)$. In particular the oberserved $\sigma$-field at time $t$ becomes $\Ob_t= \sigma(R_t, R_tY_t)$.

\section{Applications}

\subsection{Right-censoring of counting processes with time-dependent covariable}
We consider the modelling of independent counting processes $W^i=(W^i_t)$ with possibly time-dependent explanatory
variables $Z^i=(Z^i_t)$. We observe $X^i=(W^i,Z^i)$, $i=1,\ldots,n$ through a mechanism specified by $R=(R^{W^1},
\ldots,R^{W^n};R^{Z^1}, \ldots,R^{Z^n})$. The $(Z^i_t)$ are supposed to be completely observed so that
$R^{Z^i}_t=1$ for all $t$. We consider a family of probability measures $\{P^{\theta,\gamma,\psi}\}$,
where $\theta$ is the parameter of interest which parameterizes the dynamics of $(W^i_t)$ given the value of the
 explanatory variables, that is, the intensity of $(W^i_t)$ in the filtration $(\X^i_t)$ (defined as the family of $\sigma$-fields $\X^i_t=\W^i_t \vee \Z^i_t$) depends on
$\theta$ only. Let $\gamma$ parameterize the marginal law of $(Z^i_t)$. More precisely, it is assumed that the compensator of $Z^i_t$ is the same in the filtration  $(\X^i_t)$ and in the filtration $(\Z^i_t)$; in other word it is assumed that $Z$is WCLI of $W$; see Commenges and G\'egout-Petit \cite{Commenges09}. If this is not the case, then we must resort to joint modelling. Consider that assumption A1 holds; for right-censoring this means that the set of times at which observation may be stopped is denumerable; this is not a limitation in practice: for instance in an epidemiological cohort we may say that observation may be stopped each day at a fixed hour (we generally do not have a precision better than one day). If ignorability holds for $(X,R)$, we have on $\{R=r\}$,
 $\LL^{\theta, \gamma / \theta_0, \gamma_0 }_{\X^r}=
\LL^{\theta, \gamma, \psi_0 / \theta_0, \gamma_0, \psi_0}_{\Ob}$. But we have also $\LL^{\theta, \gamma /
\theta_0, \gamma_0 }_{\X^r}=\LL^{\theta/ \theta_0}_{\W^r|Z^r}\LL^{\gamma/ \gamma_0}_{Z^r}$, so that in terms of
inference about $\theta$ we only need to compute $\LL^{\theta/ \theta_0}_{\W^r|Z^r}$. The conclusion is that although we can use the conditional likelihood of $W$ given $Z$, we may still consider the
ignorability condition for $(X,R)$, that is, the response indicator process may depend on both observed $W$ and
$Z$. This discussion is related to the so-called ``reduced model'' introduce in \cite{Commengesetal07}.

We give a particular example of artificial right censoring of $0-1$ counting processes which is compatible with
ignorability (assuming to simplify that there is no other source of random censoring). Consider a study where
$n$ subjects are potentially followed-up until a time $t^*$; it is assumed that the processes $(X^i,Z^i)$ are independent and identically distributed. Suppose that at a given time $t_1$ we make an analysis
of the data. Using for instance maximum likelihood estimators in a parametric model we can construct an
estimator $\hat \theta_1$ which by definition is $\Ob_{t_1}$-measurable. Let us suppose for simplicity that
$Z_t, t>t_1$ is known at $t_1$; it is then possible to compute an estimator
of the probability that subject $i$ experiences the event before the end of the study $P_{\hat
\theta}(W^i_{t^*}=1|\W_{t_1},\Z_{t_*})$, where $\W_t$ and $\Z_t$ denote $\sigma$-fields generated by the $n$ processes up to time $t$. For reducing the cost of the study without reducing too much its power,
we may decide to follow after $t_1$ only the subjects for whom $P_{\hat \theta}(W^i_{t^*}=1|\W^i_{t_1},\Z^i_{t^*})
\ge c$, for some chosen $c$. In the TCMP this means putting $R^{W^i}_t=0$ for $t>t_1$ for those subjects with an
estimated probability below $c$. It is clear that $R$ is $\Ob_t$-predictable which  from Corollary \ref{Rprev} implies ignorability on all $r$ because A1 holds. Note that ignorability is not dependent of the good specification of
the model used at $t_1$; of course the validity of the final analysis will depend on the good specification of
the model used for it.

\subsection{Longitudinal markers (continuous state-space processes)}
\subsubsection{Missing data at fixed observation times}
This is the classical set-up of ``repeated measurements'' or ``longitudinal data'' (in a narrow sense). We
consider independent continuous state-space processes $(X^i_t)$, $i=1,\ldots,n$; for simplicity we do not consider covariates. Here $(X^i_t)$ is planned to be observed at $v^i_j$, $j=1, \ldots,
m$ but there may be missing data. The response indicator processes for
$(X^i_t)$ can be written: $R_t^{X^i}=\sum_{j=1}^m 1_{\{t=v^i_j\}}S_j^i$, where $S_j^i$ are  binary variables. The jump at $v^i_j-$ of the
compensator of the counting process associated to $R^i$ in a filtration $(\G_t)$ is $P(S_j^i=1|\G_{v^i_j-})$. CAR(DYN) can thus be expressed as:
$P(S_j^i=1|S_1^i,\ldots, S^i_{v_{j-1}^i},\X)=P(S_j^i=1|S_1^i,\ldots, S^i_{v_{j-1}^i},S_1^iX_1^i,\ldots, S^i_{v^i_{j-1}}X^i_{v^i_{j-1}})$, which implies (by taking conditional expectation) that
$$P(S^i_j=1|S^i_1,\ldots, S^i_{v^i_{j-1}},\X^i_{v^i_{j}-})=P(S^i_j=1|S^i_1,\ldots, S^i_{v^i_{j-1}},S^i_1X^i_1,\ldots, S^i_{v^i_{j-1}}X^i_{v^i_{j-1}});$$ this can intuitively be interpreted in saying that the
missing data mechanism may depend (only) on the observed $W^i$ and on the $S^i$ up just before time $v^i_j$; this case could be treated with the conventional MAR concept.

\subsubsection{Random observation times}\label{randomobservationtimes}
We may consider determining observation times in order to reduce the costs of a study or to improve the
monitoring of patients. Consider the case of a study of the evolution of the concentration of CD4+ T lymphocytes in HIV infected
patients, represented by $(X^i_t)$ for subject $i$. Following example \ref{CD4}, the observations are the CD4 counts, that is noisy measurements of $(X^i_t)$ at discrete visit times: $Y_{ij}=X^i_{V^i_{j}}+\varepsilon_{ij}$. The time of the next visit, $V^i_{j+1}$ may
depend on  CD4 count at $V^i_{j}$ according to the procedure: $V^i_{j+1}=V^i_{j}+g(X^i_{V^i_{j}})$, where $g(.)$ is a known function.. For instance we could decide
to see a patient with a delay of three months if CD4 $\ge 500$, two months if $200 \le$ CD4 $<500$ and one month
if CD4 $< 200$. It is clear that CAR(DYN) would hold in this instance. On the contrary if the visit time was decided
based in part on clinical symptoms not included in the model, or if drop-out could be due to severe clinical
events or death (related to CD4 and not included in the model), CAR(DYN) would not hold.

\subsubsection{Vertically coarsened observations}
If we are interested in the evolution of HIV-RNA, we may consider the same issues as above ($X^i_t$ now
representing HIV-RNA), with the additional complexity of a (known) detection limit $\eta)$: this produces a left-censoring of $X^i_t$:  if $R^i_t=1$ we observe  $Y^i_t=(1_{X^i_t> \eta}, \tilde X^i_t)$, where $\tilde X^i_t=\max (X^i_t,\eta)$ (see \cite{Thie}). This is a case of what we have called ``vertical coarsening'' in section 3.2.

\subsection{Joint modeling}
One of the reasons for considering a joint model is precisely to remove the bias due to so called ``informative
censoring ``. Consider as in the example of section \ref{randomobservationtimes} that we are interested in the evolution of the concentration of CD4+ T-lymphocytes
 represented by $(T^i_t)$. Only noisy discrete observation of $(T^i_t)$, the CD4 counts, are available. Assume that subjects are lost from follow-up when they develop AIDS. There is a rather strong relationship between concentration of CD4+ T-lymphocytes and the risk of developing AIDS, so that the
intensity of the counting process describing drop-out depends on it: CAR(DYN) does not
hold in a model which does not include AIDS.

Thus we may consider jointly modelling the concentration of CD4+ T-lymphocytes and AIDS and consider the process $(X^i_t)=(T^i_t,W^i_t,Z^i_t)$,
where $(W^i_t)$ is a counting process which represents AIDS ans $(Z^i_t)$ is a
multivariate process of explanatory variables (such a joint model has been developed in \cite{Guedj}). We may allow as in section \ref{randomobservationtimes} that the visit times depend on the
observed CD4 counts; we may also allow that the probability that the subject drop out after $V^i_j$ depends on
the AIDS status at this visit. The compensator of the associated counting process $N^i$ will be constant between
$V^i_j$ and $V^i_{j+1}$ and will make a jump at $V^i_{j+1}$ equal to one minus the probability of drop-out: if
the latter depends on the observed AIDS status $W^i_{V^i_j}$ only, then CAR(DYN) holds; if it depends of unobserved
status $W^i_t$ for $t\ne V_j$, then CAR(DYN) does not hold.

\section{Conclusion}
We have proposed a coarsening model for processes (TCMP) and  developed a theory of ignorability in this framework. The theory applies to general stochastic processes having discrete or continuous state-space; in particular it applies to both counting processes and diffusion processes. The framework of repeated measurements can be represented as a continuous-time continuous state-space process observed at discrete times.
Our results hold even if the observed part of the process of interest $X$ has a null probability, which is the case in the examples of the previous section. We have given a factorization condition for the likelihood which allows to get rid of the nuisance parameter and may be called weak ignorability; we can define ignorability in a strong sense, that is equality between the correct likelihood and the likelihood ignoring the m.l.i.d.  on events $\{R=r\}$ of non-null probability. This restriction comes from the fact that likelihoods, as Radon-Nikodym derivatives, are not uniquely defined.

We have already applied some of the results presented here to a multi-state model for dementia, institutionalization and death \cite{Commenges07}. This model can be represented by a three-variate counting process to which we associate the three-variate response process. In this application the observed event of $X$ is generally of null probability although it is not a singleton; institutionalization can be observed either exactly or in an interval; we showed  that the m.l.i.d. could be ignorable. The factorization theorem is also useful for defining a risk for model choice in the context of coarsened observations of multistate processes as is proposed in \cite{Commengesetal07}. Finally this scheme of observation is useful for making inference in the framework of a general dynamical model such as proposed by Commenges and G\'egout-Petit \cite{Commenges09}. In these latter papers it is made explicit that there is a true probability  $P^*$ under which the events in the universe are generated. In the present paper we have spoken in terms of model. However it must be clear that for ignorability to hold really, the conditions must hold under the true probability; in particular the condition in CAR(DYN) must be hold under $P^*$.

\appendix

\section{CAR(TCMP)-loc}\label{loc}

This appendix develops the concept of CAR(TCMP)-loc of Definition \ref{carloc}.
First note that this definition has a meaning only for $r$ such that $P(\{R=r\})>0$.
The TCMP model verifies the condition  CAR(ABS) defined by
Nielsen (2000) if $\forall (\theta, \psi)$, $P_{(\theta, \psi)}$
is CAR(ABS) i.e. the following condition is true :
\begin{eqnarray}
\mbox{for } P_{\theta} \mbox{ a.e. $x$, $x'$, for every }A \in \Ob, & &
\nonumber \\
 P_{ \psi}(A \cap D_x \cap D_{x'}|X=x)& = & P_{ \psi}(A \cap
D_x \cap D_{x'}|X=x')\label{carabs}
 \end{eqnarray}
  where $D_x=\{(r,y)\in (\Xi,\Gamma); ry=rx\}\in \Ob$. Note moreover that as
  pointed by Gill \& all \cite{Gill}, if $P_{(\theta_0, \psi_0)}$ is
  CAR(ABS) and the model is CAR(REL), then the model is CAR(ABS).

\begin{thm} If the TCMP is CAR(ABS), then CAR(TCMP)-loc holds on all $r$.
\end{thm}

\noindent {\bf {\em Proof.}} Assume the model is CAR(ABS) and remark that $\{R=r\}
\cap D_x=\{R=r\} \cap D_{rx}$ and that $\forall (\theta, \psi)$,
$P_{(\theta, \psi)}[{D_x}|X=x]=1$, then $\forall (\theta, \psi)$
for $P_{\theta}$ a.e. $x$
 and every $\Ob$-measurable $Z$, we get
\begin{eqnarray*}
\EZZZ[1_{\{R=r\}}Z|X=x]& = & \EZZZ[1_{\{R=r\}}Z 1_{D_x}|X=x]\\
& = & \EZZZ[1_{\{R=r\}}Z |X=rx]
\end{eqnarray*}
From this equality, we deduce that $\forall (\theta, \psi)$ for
 every $Z$, $\Ob$-measurable, $\EZZZ[1_{\{R=r\}}Z|\X] = \EZZZ[1_{\{R=r\}}Z
|\X^r]$. It follows that for every $A \in \R$:
$\EZZ[1_{\{R=r\}}1_A {\cal L}^{(\theta, \psi) /(\theta_0,
\psi_0)}_{\R|{\cal X}}|\X] = \EZZZ[1_{\{R=r\}}1_A {\cal
L}^{(\theta, \psi) /(\theta_0, \psi_0)}_{\R|{\cal X}^r}|\X]$ which
implies CAR(TCMP)-loc. \hfill \ding{113}

\noindent {\bf {\em Proof of Theorem \ref{loc}.}} By the iterated decomposition formula (property v), we have ${\cal L}^{(\theta,
\psi) /(\theta_0, \psi_0)}_{\X,\R}={\cal L}_{\X,\X^r,\R}={\cal L}^{(\theta, \psi) /(\theta_0,
\psi_0)}_{\X|\X^r,\R}{\cal L}_{\R|\X^r}{\cal L}_{\X^r}$. We have ${\cal L}^{(\theta, \psi) /(\theta_0,
\psi_0)}_{\Ob}={\rm E}[{\cal L}^{(\theta, \psi) /(\theta_0, \psi_0)}_{\X,\R}|\Ob]={\rm E}[{\cal L}^{(\theta,
\psi) /(\theta_0, \psi_0)}_{\X|\X^r,\R}{\cal L}_{\R|\X^r}{\cal L}_{\X^r}|\Ob]$. On $\{R=r\}$, $\Ob= \X^r \vee
\R$ and by  Lemma \ref{Kallenberg} we have:
$${\cal L}^{(\theta, \psi) /(\theta_0, \psi_0)}_{\Ob}=
{\cal L}^{\theta/\theta_0}_{\X^r}{\cal L}^{(\theta, \psi)
/(\theta_0, \psi_0)}_{\R|\X^r}{\rm E}_{(\theta_0, \psi_0)} [{\cal
L}^{(\theta, \psi) /(\theta_0, \psi_0)}_{\X|\X^r,\R}|\X^r\vee
\R]={\cal L}^{\theta/\theta_0}_{\X^r}{\cal L}^{(\theta, \psi)
/(\theta_0, \psi_0)}_{\R|\X^r}.$$ If CAR(TCMP)-loc holds we have
on $\{R=r\}$, ${\cal L}^{(\theta, \psi) /(\theta_0,
\psi_0)}_{\R|{\cal X}^r}={\cal L}^{(\theta, \psi) /(\theta_0,
\psi_0)}_{\R|{\cal X}}$; we have ${\cal L}^{(\theta, \psi)
/(\theta_0, \psi_0)}_{\R|{\cal X}}={\cal L}^{\psi
/\psi_0}_{\R|{\cal X}}$. Thus with both conditions we have ${\cal
L}^{(\theta, \psi) /(\theta_0, \psi_0)}_{\Ob}={\cal
L}^{\theta/\theta_0}_{\X^r}{\cal L}^{\psi/\psi_0}_{\R|\X},$ a.s.
on ${\{R=r\}}$
 and this implies ignorability.\hfill \ding{113}
\vspace{6mm}

\section{Proof of Theorem \ref{idpcensor}}

 The likelihood of the counting process $(X, N)$ can be written using Jacod's formula as
$\LL_{\X,\N}=\phi(\Lambda^{\FF,X},\Lambda^{\FF,X}_0,X)\phi(\Lambda^{\FF,N}, \Lambda^{\FF,N}_0,N)$. As we have done above we can also write it:
$\LL_{\X,\N}=\LX \phi(\Lambda^{\FF^*,N}, \Lambda^{\FF^*,N}_0,N)$. Noting that $\LX = \phi(\Lambda^{\X,X},\Lambda^{\X,X}_0,X)$ and equating the two
representations we have:
$$\phi(\Lambda^{\FF,X}, \Lambda^{\FF,X}_0,X)\phi(\Lambda^{\FF,N}, \Lambda^{\FF,N}_0,N)=\phi(\Lambda^{\X,X}, \Lambda^{\X,X}_0, X) \phi(\Lambda^{\FF^*,N} \Lambda^{\FF^*,N}_0,N).$$
CAR(DYN) says that $\Lambda^{\FF^*,N}=\Lambda^{\Ob,N}$ and in this
right-censoring case we have $\Lambda^{\Ob,N}=\Lambda^{\FF,N}$
(this is because $\Ob_t=\FF_t$ on $\{t\le C\}$ and
$\Lambda^{N}_{C+u}=\Lambda^{N}_{C}$ whatever the filtration). So
if CAR(DYN) holds, the above equation yields
$\phi(\Lambda^{\FF,X}, \Lambda^{\FF,X}_0,
X)=\phi(\Lambda^{\X,X},\Lambda^{\X,X}_0,X)$. This  must be true
almost surely, for all $(\theta, \psi)$, and moreover, we still
have this equality if we stop the observation at time $t$ or at a
$(\X_t)$-stopping time $T$. All that we have to prove is that this
implies $\Lambda^{\FF,X}=\Lambda^{\X,X}$, which is independent
censoring.

Let us begin with $X$ a $0-1$ counting process and denote its jump time $T$. If we stop observation at $t$, we have on $\{T>t\}$, $\frac{\exp \Lambda^{\FF,X}_t}{\exp \Lambda^{\FF,X}_{0t}}=\frac{\exp \Lambda^{\X,X}_t}{\exp \Lambda^{\X,X}_{0t}}$; because of left-continuity, we have also the equality of $t=T$ and because the intensity is equal to zero after $T$, the equality holds for all $t$ almost surely. Taking log and differentiating we obtain:
\begin{equation}  \lambda^{\FF,X}_t - \lambda^{\FF,X}_{0t}=\lambda^{\X,X}_t - \lambda^{\X,X}_{0t}.\label{differ}\end{equation}
The likelihood has a limit when  $t \rightarrow \infty$ and at the limit we have $\frac{\lambda^{\FF,X}_T\exp \Lambda^{\FF,X}_T}{\lambda^{\FF,X}_{0T}\exp \Lambda^{\FF,X}_{0T}}=\frac{\lambda^{\X,X}_T\exp \Lambda^{\X,X}_T}{\lambda^{\X,X}_{0T}\exp \Lambda^{\X,X}_{0T}}$ from which we successively deduce
$\frac{\lambda^{\FF,X}_T}{\lambda^{\FF,X}_{0T}}=\frac{\lambda^{\X,X}_T}{\lambda^{\X,X}_{0T}}$ and
$\frac{\lambda^{\FF,X}_t}{\lambda^{\FF,X}_{0t}}=\frac{\lambda^{\X,X}_t}{\lambda^{\X,X}_{0t}}$, almost surely for all $t$ on the support of the distribution of $T$. Combining this result with (\ref{differ}), we obtain $\lambda^{\FF,X}_t =\lambda^{\X,X}_t $ $a.s.$, which for c\`ag processes implies indistinguishability of the intensities and of the cumulative intensities. If the process may have several jumps $T_1, T_2, \ldots$, we first prove by the same reasoning that we have equality of the intensities on  $\{t\le T_1\}$, then using this result and again the same reasoning we have equality on $]T_1, T_2]$ and so on.
All this reasoning is symmetrical so we can also prove that independent censoring implies CAR(DYN). \hfill \ding{113}
\vspace{6mm}

\section{Extension to left-continuous $R$}

In some situations it is natural to consider response indicator
processes which are left-continuous; case II right-censoring is an
example. In that case it is not possible to directly associate to
$R$ a counting process and hence to apply CAR(DYN). In order to
extend the application of CAR(DYN) to such processes, we will
consider them as limits of right-continuous processes. Let us
still consider the univariate case. Consider for instance the case
where $R_0=1$. $(R_t)$ may be left-continuous at jumps of odd
ranks: $V_{2j+1}, j\ge 0$; the process can be written
$R_t=\sum_{j\ge 0}1_{[V_{2j},V_{2j+1}]}$. Consider the sequence of
processes ${R^n_t}=(\sum_{j\ge 0}1_{[V_{2j}\le t <V^n_{2j+1}[})
\wedge 1$, defined by: $R^n_0=1$, $V^n_{2j+1}=V_{2j+1}+1/n$. The
limit of $(R^n_t)$ is $(R_t)$.
\begin{thm} \label{Rcag}
Consider a process $R=(R_t)$ which is right-continuous at upward
jumps and may left-continuous at downward jumps. Consider a
sequence of right-continuous processes $R^n=(R^n_t)$ constructed
as above; if each $(X,R^n)$ satisfies  CAR(TCMP)-loc on $r$
then ignorability holds for $(X,R)$ on $r$.
\end{thm}

\noindent {\bf {\em Proof.}} We note  $\Ob^n$, the observed $\sigma$-field
associated to ${R^n}$. If $(X,R^n)$ satisfies CAR(TCMP)-loc on
$r$ then (see the proof of Theorem (\ref{loc})), we have for
all $n$: on $\{R=r\}$,
$\LL^{(\theta,\psi)/(\theta_0,\psi_0)}_{\Ob^n}=
\LL^{(\theta,\psi)/(\theta_0,\psi_0)}_{\R^n|\X}\LL^{\theta/\theta_0}_{\X^{r^n}}$.
 $\Ob^n$ is
larger than $\Ob$: $\Ob \subset \Ob^n$ and it is clear that
$\Ob^n$ is a decreasing sequence of $\sigma$-fields: $\Ob= \cap_n
\Ob^n= \Ob^{\infty}$. By the Downward Levy Theorem \cite{Williams} we have: $\LL^{(\theta,\psi)/(\theta_0,\psi_0)}_{\Ob^n}={\rm
E}_{(\theta_0,\psi_0)}(\LL^{(\theta,\psi)/(\theta_0,\psi_0)}_{\FF}|\Ob^n)
\rightarrow {\rm
E}_{(\theta_0,\psi_0)}(\LL^{(\theta,\psi)/(\theta_0,\psi_0)}_{\FF}|\Ob)=
\LL^{(\theta,\psi)/(\theta_0,\psi_0)}_{\Ob} $ $a.s.$. Using again
the Downward Levy Theorem we get $\LL^{\theta/\theta_0}_{\X^{r^n}}
\rightarrow \LL^{\theta/\theta_0}_{\X^{r}}$ a.s. Moreover note
that $\R^n=\R$ because the process $R^n$ is deterministically
defined from $R$: this implies that
$\LL^{(\theta,\psi)/(\theta_0,\psi_0)}_{\R^n|\X}=\LL^{\psi/\psi_0}_{\R|\X}$.
At the limit we have thus:
$\LL^{(\theta,\psi)/(\theta_0,\psi_0)}_{\Ob}=
\LL^{\psi/\psi_0}_{\R|\X}\LL^{\theta/\theta_0}_{\X^{r}}$ which
concludes the proof. A similar result could be obtained for upward
jumps. \hfill \ding{113}

 As an example consider the case of Type II right-censoring
where we have independent $0-1$ counting processes $(X^i_t)$, $i=1,\ldots,n$ and observation is stopped just
after observing the $d^{th}$ event. Thus the response indicator process $R$ is not independent on the
multivariate process $X$. In fact we have $R_t=1_{\{\bar X_{t}\le d\}}$, where $\bar X_t=\sum_i X^i_t$. This is
a case of a left-continuous process which has only one downward jump. Since $R$ is $\X$-measurable $\LRX=1$ (by
property iv) and is thus obviously $\Ob$-measurable which is CAR(TCMP). Consider now a slightly more
sophisticated model which we call randomized Type II censoring in which we may stop observation after each event
with a given probability depending of what have been observed. For instance let $(T_1, T_2,...,T_n)$ the times
of occurrence of the first, second,..., events, and let the probability of stopping observation just after $T_j$
(conditional on having observed $X$ until $T_j$) be $\frac{j-1}{j}$, $j=1,\ldots,n$; let $C$ be the jump time of
$R$ ($C=T_j$, for some $j$).
 We consider $R$ as the limit of the sequence of right-continuous processes $R^n$ such that $R^n_t=1_{\{t\ge C+1/n\}}$. We
can easily verify that these observation processes satisfy CAR(DYN) (because future values of $X$ are not used
for defining the probability of stopping observation), and thus CAR(TCMP) by Theorem 5; thus $R$ itself
satisfies CAR(TCMP) by Theorem \ref{Rcag}.
\vspace{6mm}

\end{document}